\documentclass[conference]{IEEEtran}
\IEEEoverridecommandlockouts

\usepackage{cite}
\usepackage{amsmath,amssymb,amsfonts}
\usepackage{algorithmic}
\usepackage{graphicx}
\usepackage{textcomp}
\usepackage{xcolor}
\usepackage{booktabs}
\usepackage{url}
\usepackage{hyperref}
\usepackage{enumitem}

\hypersetup{
    colorlinks=false,
    linkcolor=blue,
    citecolor=blue,
    urlcolor=blue
}

\def\BibTeX{{\rm B\kern-.05em{\sc i\kern-.025em b}\kern-.08em
    T\kern-.1667em\lower.7ex\hbox{E}\kern-.125emX}}

\begin{document}

\title{Mapping Austria's Natural Gas and Hydrogen Infrastructure Plans}

\author{

\IEEEauthorblockN{
Marco Quantschnig,
Yannick Werner,
Thomas Klatzer and
Sonja Wogrin
}

\IEEEauthorblockA{
Institute of Electricity Economics and Energy Innovation
}

\IEEEauthorblockA{
Research Center for Energy Economics and Energy Analytics \\ Graz University of Technology \\ Graz, Austria
}

\IEEEauthorblockA{
\{marco.quantschnig, yannick.werner, thomas.klatzer, wogrin\}@tugraz.at
}

}

\maketitle

\begin{abstract}
This paper presents a comprehensive, spatially disaggregated dataset of Austria's natural gas and hydrogen infrastructure towards 2040. The dataset covers the complete gas transmission and distribution networks down to the medium-pressure level and integrates hydrogen expansion plans from the Austrian Gas Grid Management. Transmission infrastructure is reconstructed from ENTSOG maps, converted into a topologically consistent graph representation, and enriched with technical attributes through automated spatial matching with open-source datasets such as OpenStreetMap and Global Energy Monitor. Distribution networks and infrastructure modifications are implemented using QGas, a newly developed GIS-based tool for graph-based infrastructure manipulation. To enable forward-looking energy system analyses, the dataset explicitly represents the stage-wise transition from natural gas to hydrogen infrastructure within one single dataset. Repurposed and newly constructed hydrogen pipelines are integrated within a unified network topology using node-splitting and time-dependent connector elements, enabling consistent modeling of parallel natural gas and hydrogen operation over time. The resulting dataset provides a detailed representation of Austria's gas and hydrogen infrastructure, including 586 natural gas pipeline segments (5{,}000 km), 113~repurposed segments (1{,}250 km), and 39 newly constructed hydrogen segments (820 km), connecting 720 nodes. Moreover, it includes a comprehensive set of gas demands, biogas production facilities, storage units, electrolyzers, and compressor elements, making it directly applicable for energy system optimization models.
\end{abstract}


\section{Introduction}

Achieving climate neutrality requires a fundamental transformation of energy systems, in which large-scale electrification will be complemented by renewable gases to decarbonize hard-to-abate processes and provide long-duration energy storage~\cite{kharel2018}.
Energy system optimization models (ESOMs) are widely used tools to support the planning and decision-making processes associated with this transformation. To this end, ESOMs require spatially resolved representations of energy infrastructures, and in particular network topologies, as this is decisive for the validity of model results at national and regional scales~\cite{zwickl2024,aryanpur2021}. This is highlighted in recent integrated power and hydrogen system optimization studies~\cite{neumann2023,reuss2019}, which show that model outcomes and investment decisions, e.g., regarding the repurposing of natural gas pipelines for hydrogen transmission and the expansion of new hydrogen pipelines, are highly sensitive to the spatial detail of gas network representations.

While publicly accessible information on gas network infrastructure exists, it is typically not provided in a form suitable for direct use in ESOMs. Data is often available in heterogeneous formats and at varying levels of detail, which limits their applicability. For example, ENTSOG publishes PDF-based infrastructure maps~\cite{entsog} that illustrate pipeline routes and major facilities but are not directly usable for quantitative analyses. To address these limitations, open datasets such as SciGRID\_gas~\cite{pluta2022} convert these maps into topologically consistent representations of the European gas infrastructure. However, this dataset is primarily designed for case studies at the European scale and therefore lacks the spatial resolution and adaptability needed for detailed national or regional analyses.

For Austrian case studies, these limitations are particularly consequential. The Austrian gas network features a dense distribution system and multiple parallel pipeline loops~\cite{tag2025,gasconnect2025}, which are critical for sequential repurposing strategies to maintain natural gas supply during the transition. National expansion plans published by the Austrian Gas Grid Management (AGGM) reflect this approach~\cite{aggm_h2_roadmap,aggm_knep}. Existing Austrian datasets, such as Zwickl-Bernhard et al.~\cite{zwickl2024,zwickl2023zenodo}, do not represent parallel pipelines as separate loops (e.g., the Trans Austria Gasleitung), which limits their applicability for hydrogen system optimization studies. To the best of our knowledge, no open-source dataset currently exists that provides the gas infrastructure representation to integrate and analyze Austria's hydrogen expansion plans in a highly spatially resolved and temporally consistent graph format suitable for ESOMs.

To address this gap, this paper presents an open-source dataset of Austria's natural gas and hydrogen infrastructure intended to support quantitative energy system analyses. The dataset provides a consistent representation of gas infrastructure and its potential evolution over time, enabling long-term and long-horizon case studies of energy system transformation.

Specifically, the dataset presented in this paper makes three key contributions:

\begin{enumerate}[leftmargin=1.5em]
    \item It provides an open-source and highly detailed representation of Austria's natural gas and hydrogen infrastructure, ready to be used in ESOMs. The dataset features a fully resolved network topology that explicitly represents individual pipeline strands and distribution networks down to the medium-pressure level.
    \item It provides a geo-referenced representation based on standardized geodata elements (GeoJSON), enabling straightforward modification, extension, and adaptation of the dataset for the analysis of alternative expansion plans and diverse case studies. This flexibility is further enhanced by the forthcoming open-source release of the QGas tool~\cite{quantschnig2025thesis}, which allows users to efficiently manipulate and extend the network topology within a GIS-based environment.
    \item It provides a time-dependent network representation that enables a staged transition from natural gas to hydrogen infrastructure within a single unified dataset, allowing consistent long-term energy system analyses without relying on separate snapshot representations.
\end{enumerate}

The remainder of this paper is structured as follows. Section~\ref{sec:Method} introduces the general methodology and workflow for constructing and modifying the gas infrastructure dataset. Section~\ref{sec:StatusQuo} then describes the construction of the status-quo natural gas network topology and the assignment of technical parameters. Building on this topology, Section~\ref{sec:Implementation} details the integration of Austria's hydrogen expansion plans, including a time-dependent repurposing and network transition mechanism, resulting in a coherent, integrated natural gas and hydrogen infrastructure dataset. Finally, Section~\ref{sec:Validation} validates the dataset and concludes the paper.

\section{Methodology and Workflow}
\label{sec:Method}

\begin{figure*}[t]
    \centering
    \includegraphics[width=0.9\textwidth]{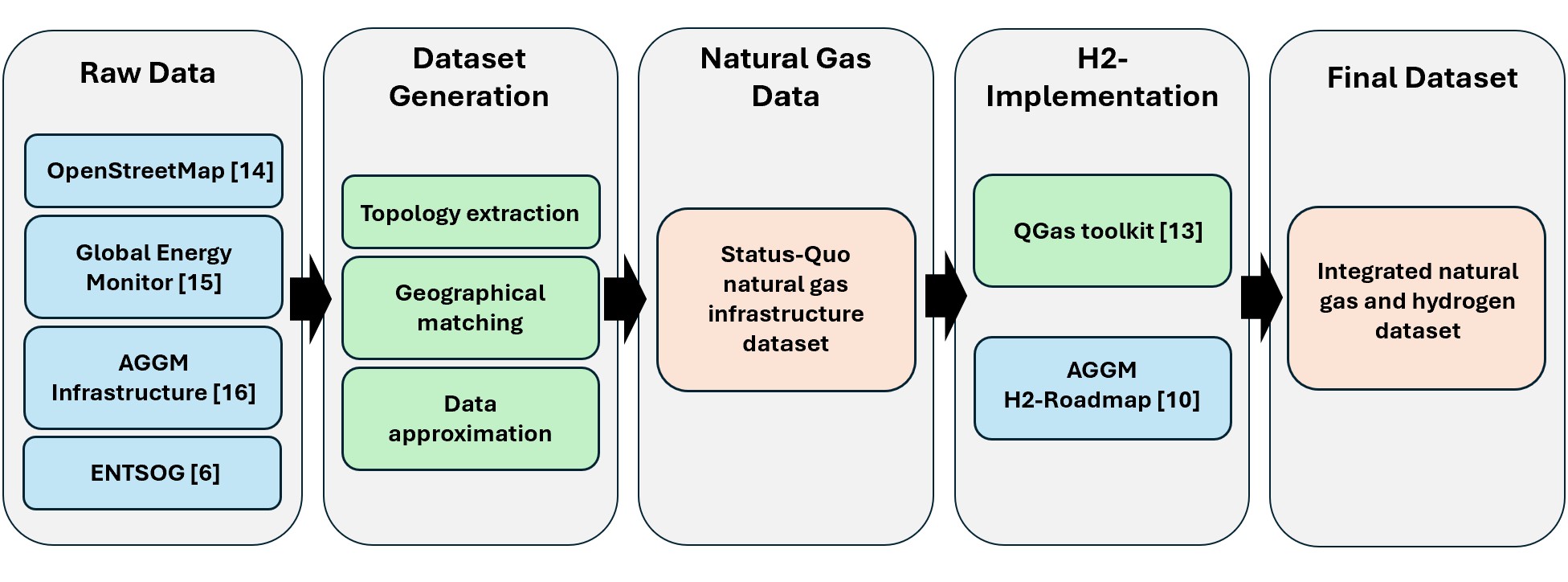}
    \caption{Schematic workflow including data sources (blue), processing steps (green), and datasets (orange).}
    \label{fig:workflow}
\end{figure*}

The process for creating the integrated natural gas and hydrogen infrastructure dataset follows the structured workflow illustrated in Fig.~\ref{fig:workflow}.

First, Austria's natural gas transmission network was reconstructed from publicly available infrastructure information. The available data were processed and harmonized to obtain a consistent representation of the transmission network suitable for quantitative energy system analyses. This step establishes the structural basis for subsequent attribute assignment, network extension, and integration of additional infrastructure layers.

Second, the transmission network was enriched with technical attributes by spatially matching it with open-source datasets. This step enables a reproducible assignment of pipeline characteristics (e.g., name, diameter, pressure) and supports future dataset updates as underlying data sources evolve. The detailed attribute matching procedure is described in Section~\ref{sec:StatusQuo}.

Third, distribution-level infrastructure was added to increase spatial resolution. High- and medium-pressure distribution pipelines were integrated based on publicly available Austrian infrastructure information using a GIS-based workflow. Due to limited public data availability at the distribution level, consistent approximations were applied. Section~\ref{sec:StatusQuo} details the integration and parameterization of distribution networks.

Finally, Austria's hydrogen infrastructure expansion plans outlined in the Hydrogen Roadmap~\cite{aggm_h2_roadmap} were incorporated in the network topology while explicitly accounting for the time-dependent network evolution and separation of emerging hydrogen and natural gas systems. The detailed implementation of hydrogen infrastructure, including the temporal transition mechanism, is described in Section~\ref{sec:Implementation}.

GIS-based infrastructure integration and data modifications were carried out using the novel QGas tool~\cite{quantschnig2025thesis}, which is employed for the visualizations in the following sections.

\section{Generation of the Status-Quo Gas Infrastructure Dataset}
\label{sec:StatusQuo}

The status-quo gas infrastructure dataset was created through a multi-step process relying exclusively on publicly available data sources. Austria's natural gas transmission network was derived from the ENTSOG Transparency Platform~\cite{entsog}, which provides infrastructure information in the form of static PDF maps illustrating pipeline routes, border points, and major facilities. These maps were first georeferenced to establish spatial consistency.

Pipeline routes were subsequently manually digitized and converted into a graph-based network topology consisting of nodes and edges. This procedure ensured topological consistency and provided a basis for subsequent network extensions and modifications. The resulting transmission network topology was then enriched with technical attributes.

Pipeline diameters, operating pressure ranges, and additional technical parameters were assigned by geographically matching the topology with open-source datasets from OpenStreetMap~\cite{osm} and the Global Energy Monitor~\cite{gem}. The automated matching procedure uses buffer-based spatial proximity by identifying corresponding infrastructure elements based on the maximum overlap length between the pipeline geometry and candidate pipelines within predefined buffer zones~\cite{quantschnig2025thesis}. To verify the accuracy of the automated matching, the assigned attributes were manually reviewed and validated using additional publicly available sources where necessary~\cite{tag2025,gasconnect2025}. Importantly, the matching workflow enables automated re-assignment of attributes using updated versions of the underlying open-source datasets, facilitating regular dataset updates and supporting future extensions to larger geographic scopes.

To further increase spatial resolution, distribution networks were added. Distribution Level~1 (high-pressure) and Level~2 (medium-pressure) pipelines were integrated based on the Austrian Infrastructure Plan published by AGGM~\cite{aggm_infra} by digitizing pipeline elements. A dedicated workflow was used to support network modification and extension. This workflow was implemented in the QGas tool~\cite{quantschnig2025thesis}, which enables the tracing and digitization of pipeline routes directly from projected infrastructure plans within the environment of the tool. Due to the limited availability of public technical data for distribution-level infrastructure, the following conservative and technically feasible approximations were applied.

Distribution Level~1 pipelines directly connected to the transmission network were represented using typical diameter ranges between 500 and 600~mm. Level~1 pipelines located further downstream, connected only to other Level~1 or Level~2 pipelines, were modeled using a representative diameter of 300 to 400~mm. Operating pressures for Level~1 pipelines were assumed to lie within a range of 20~bar to 70~bar, reflecting commonly reported high-pressure distribution conditions~\cite{econtrol,dvgw2021}.

Distribution Level~2 pipelines were represented using a diameter of 100 to 200~mm and operating pressures between 6~bar and 16~bar, corresponding to typical medium-pressure distribution networks. These values are based on commonly reported ranges for high- and medium-pressure gas infrastructure~\cite{dvgw2021,gwg2011} and were applied in a transparent and consistent manner across the dataset. All assumed parameters can be modified on an element-wise basis as more detailed information becomes available.

\begin{figure}[t]
    \centering
    \includegraphics[width=0.49\textwidth]{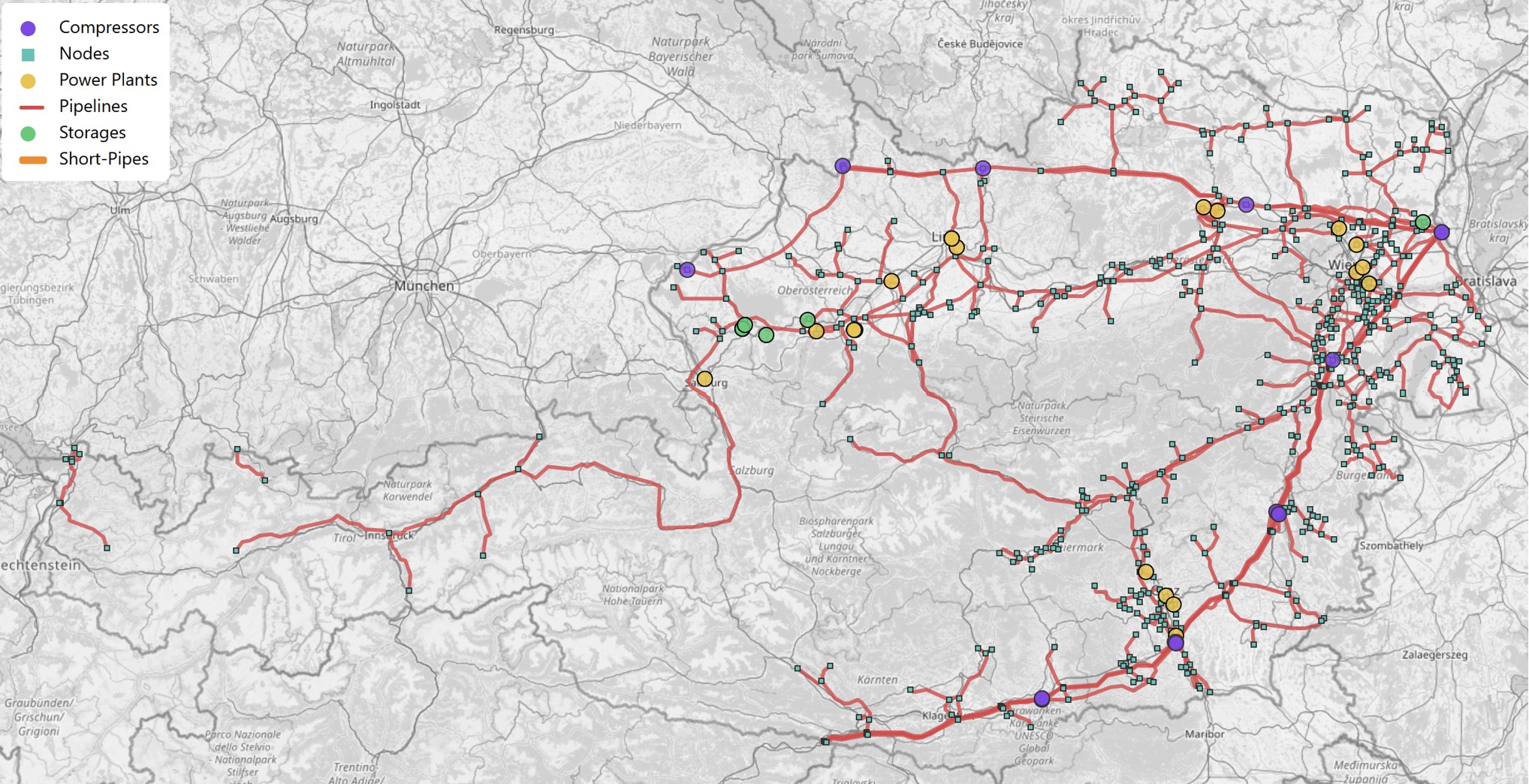}
    \caption{Status-quo natural gas infrastructure dataset for Austria.}
    \label{fig:statusquo}
\end{figure}

Fig.~\ref{fig:statusquo} illustrates the resulting status-quo gas infrastructure dataset, including transmission and distribution pipelines, network nodes, and other gas infrastructure.

\section{Implementation of Hydrogen Expansion Plans}
\label{sec:Implementation}

Austria's hydrogen infrastructure expansion plans are outlined in the AGGM Hydrogen Roadmap~\cite{aggm_h2_roadmap}, which describes a staged transition of the existing natural gas infrastructure towards hydrogen. The roadmap distinguishes four development stages ranging from 2027 to 2040 and differentiates between repurposed natural gas pipelines and newly constructed hydrogen pipelines.

\begin{figure}[t]
    \centering
    \includegraphics[width=0.49\textwidth]{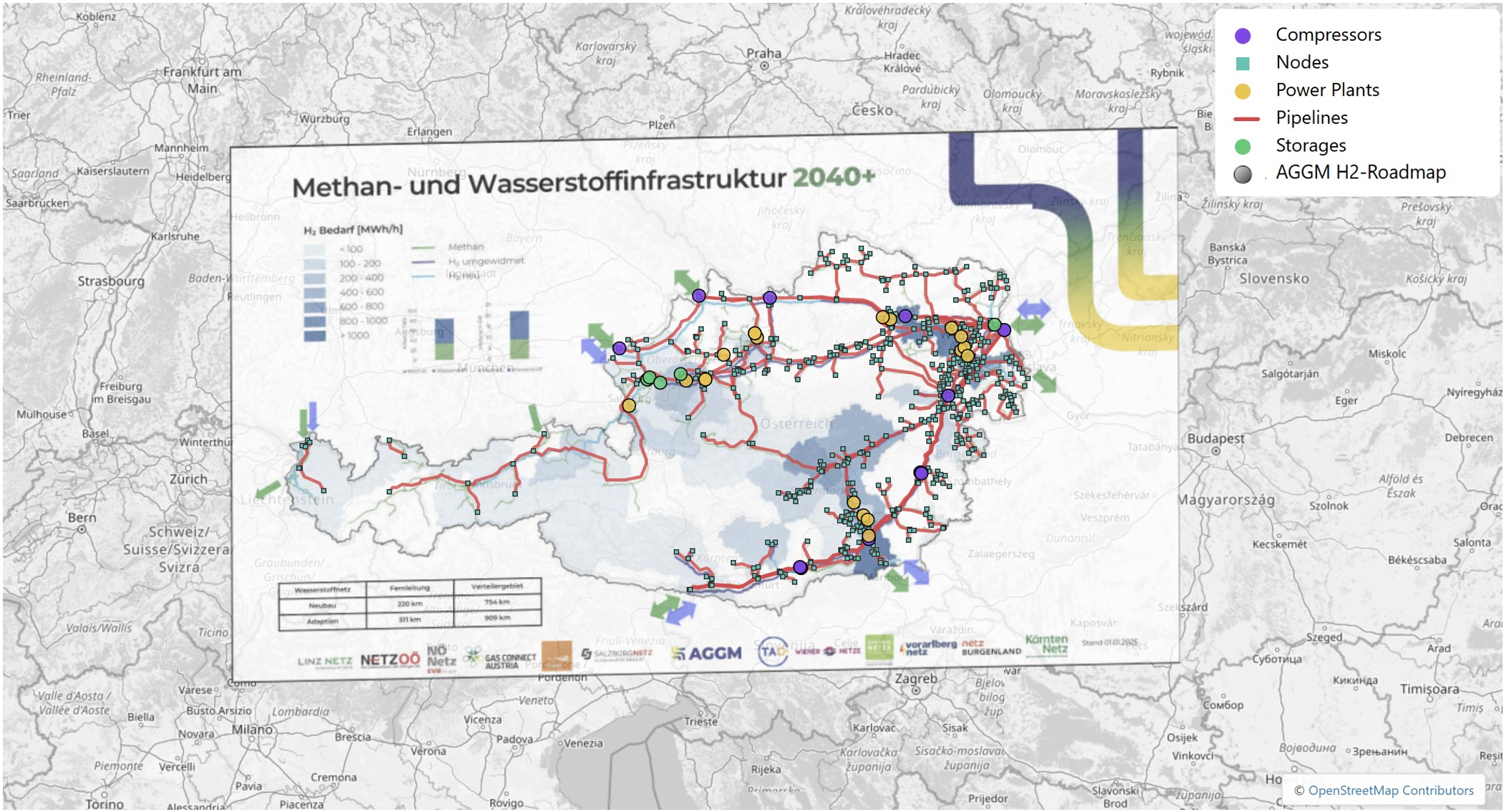}
    \caption{Projection of the Hydrogen Roadmap~\cite{aggm_h2_roadmap} in QGas.}
    \label{fig:h2roadmap}
\end{figure}

The roadmap is published as a set of graphical infrastructure maps. These maps were integrated into the QGas tool by loading them as semi-transparent background layers. The tool provides an image georeferencing functionality that allows the alignment of graphical plans with the underlying spatial reference by selecting clearly identifiable reference points in both the image and the background map. After georeferencing, the roadmap image is displayed as a separate layer, enabling the digitization of new infrastructure elements by tracing the georeferenced pipeline routes. Fig.~\ref{fig:h2roadmap} shows the projection of the georeferenced Hydrogen Roadmap stage~2040 imported into the QGas tool.

Following the georeferencing step, two dedicated pipeline sub-layers were created: one for repurposed pipelines and one for newly constructed hydrogen pipelines. Pipeline elements from the natural gas network that are designated for hydrogen conversion were transferred from the main natural gas pipeline layer into the repurposed hydrogen sub-layer. This resulted in a separate hydrogen pipeline layer with an independent and fully customizable attribute set tailored to hydrogen-specific modeling requirements. Newly planned hydrogen pipelines were digitized directly within the second sub-layer. In this stepwise manner, both repurposing plans and new infrastructure expansion plans were systematically integrated into the dataset.

\begin{figure}[t]
    \centering
    \includegraphics[width=0.49\textwidth]{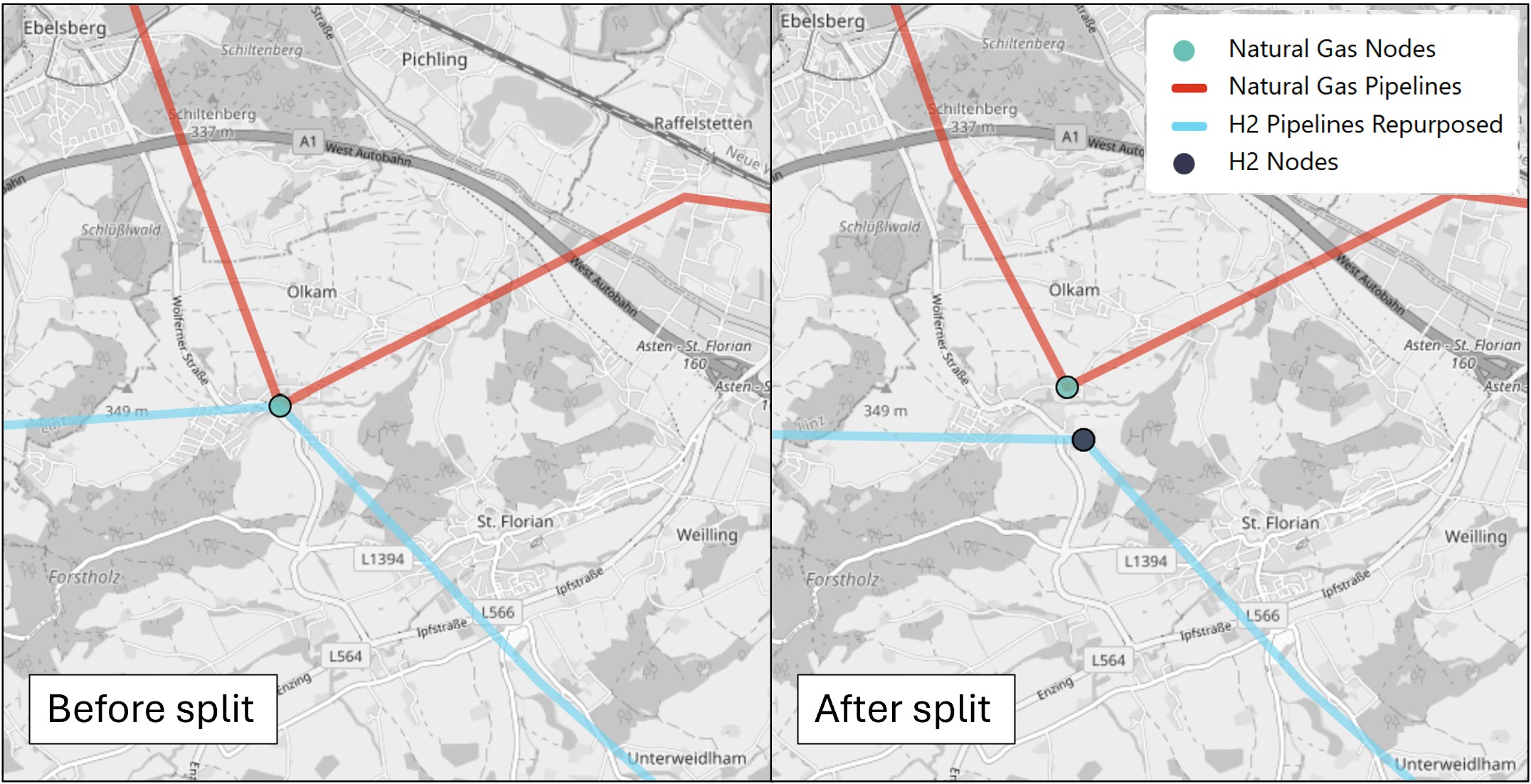}
    \caption{Splitting of existing nodes into hydrogen- and natural-gas-specific nodes.}
    \label{fig:nodesplit}
\end{figure}

Applying the methodology described above resulted in nodes at which both natural gas and hydrogen pipelines were connected. To ensure a physically consistent network representation and to prevent any artificial coupling between the two systems, these nodes were explicitly decoupled by splitting them into two sub-nodes: one assigned exclusively to the natural gas network and one to the hydrogen network. This node-splitting procedure enforces a strict separation of the two infrastructures, eliminates unintended cross-network flows in subsequent modeling applications, and yields two physically independent yet topologically consistent networks. The node-splitting process is illustrated in Fig.~\ref{fig:nodesplit}.

Each hydrogen pipeline was assigned either a repurposing year or a commissioning year, indicating the point in time at which the respective pipeline becomes operational for hydrogen. However, adjacent pipeline segments may be repurposed or commissioned at different points in time according to the Hydrogen Roadmap. In this paper we tackle this by implementing one single time-dependent network topology to avoid separate datasets which only represent a snapshot for each timestep.

\begin{figure}[t]
    \centering
    \includegraphics[width=0.49\textwidth]{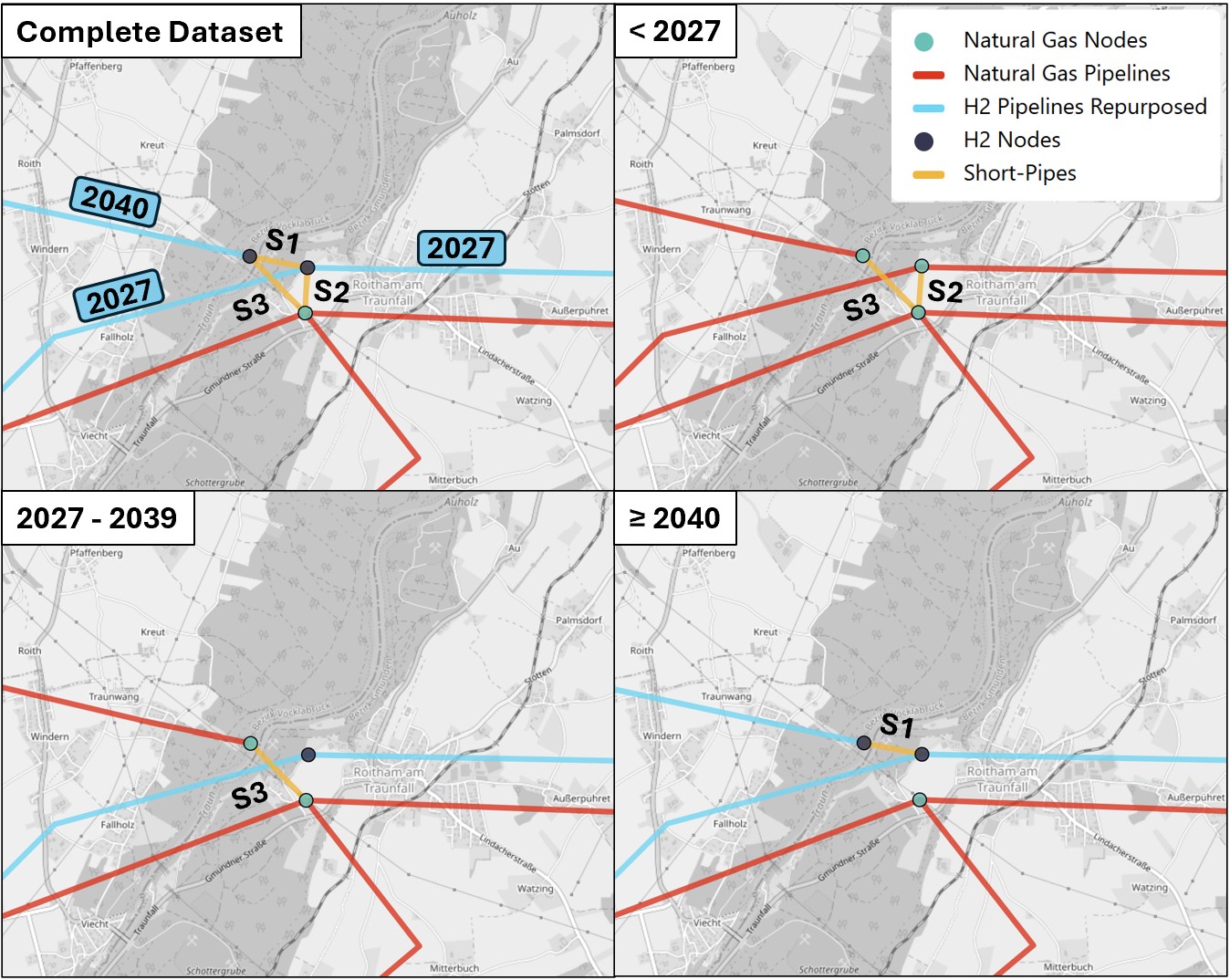}
    \caption{Visualization of the network transition.}
    \label{fig:transition}
\end{figure}

This was achieved by combining explicit repurposing and commissioning years with the introduction of short-pipe elements. These short-pipe elements represent non-limiting, lossless connectors between natural gas and hydrogen nodes and are activated or deactivated at predefined timesteps. This mechanism enables individual pipeline segments to transition from the natural gas network to the hydrogen network within a single dataset, resulting in a time-dependent change of the underlying network topology while maintaining a technically consistent and fully decoupled representation at each timestep. The underlying principle of this time-dependent transition is shown in a concrete example in Fig.~\ref{fig:transition}. Here, short-pipe S1 is activated in 2040, while S2 and S3 are deactivated in 2027 and 2040, respectively, demonstrating how substantial topological changes emerge over time while maintaining a consistent modeling framework.

Finally, the Hydrogen Roadmap also provides information on expected regional hydrogen demands through color-coded background shading of administrative districts. Based on this information, expected hydrogen demands were derived for each NUTS-3 region and assigned to demand points. These demand points contain the aggregated regional demands and were spatially positioned at the centroids of the respective NUTS-3 regions.

\section{Validation and Conclusions}
\label{sec:Validation}

The dataset was carefully validated at each timestep to ensure technical consistency: no unintended connections between hydrogen and natural gas networks exist, and the formation of isolated subnetworks is prevented, except for the import-fed Tyrol and Vorarlberg region. According to the AGGM Hydrogen Roadmap, approximately 1{,}420~km of pipelines are planned for repurposing and 940~km of new hydrogen pipelines are planned. In comparison, the presented dataset contains 1{,}253~km of repurposed pipelines and 814~km of newly constructed hydrogen pipelines. Minor deviations are due to straight-line approximations of otherwise meandering pipeline routes. In addition, the dataset includes extensive data on biogas plants, existing electrolyzers in Austria, gas storage facilities, gas-fired power plants, and compressor stations, enabling an integrated analysis of the hydrogen and natural gas infrastructure.

\begin{figure*}[t]
    \centering
    \includegraphics[width=0.92\textwidth]{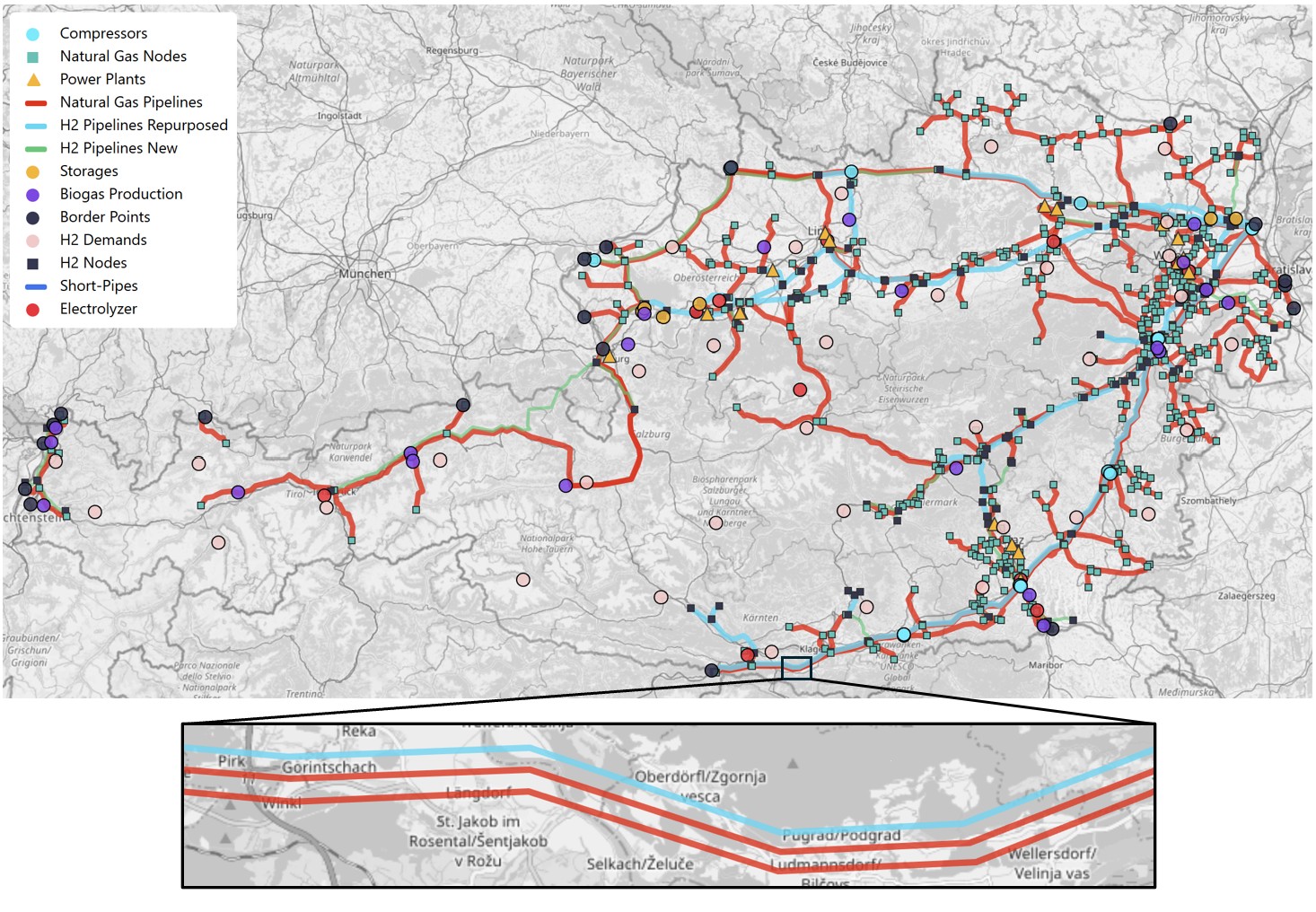}
    \caption{Integrated dataset of Austria's natural gas and hydrogen infrastructure.}
    \label{fig:integrated}
\end{figure*}

The final integrated dataset is visualized in Fig.~\ref{fig:integrated}, highlighting the importance of a dedicated representation of looped pipelines for a staged transition from natural gas to hydrogen.

\begin{figure}[t]
    \centering
    \includegraphics[width=0.49\textwidth]{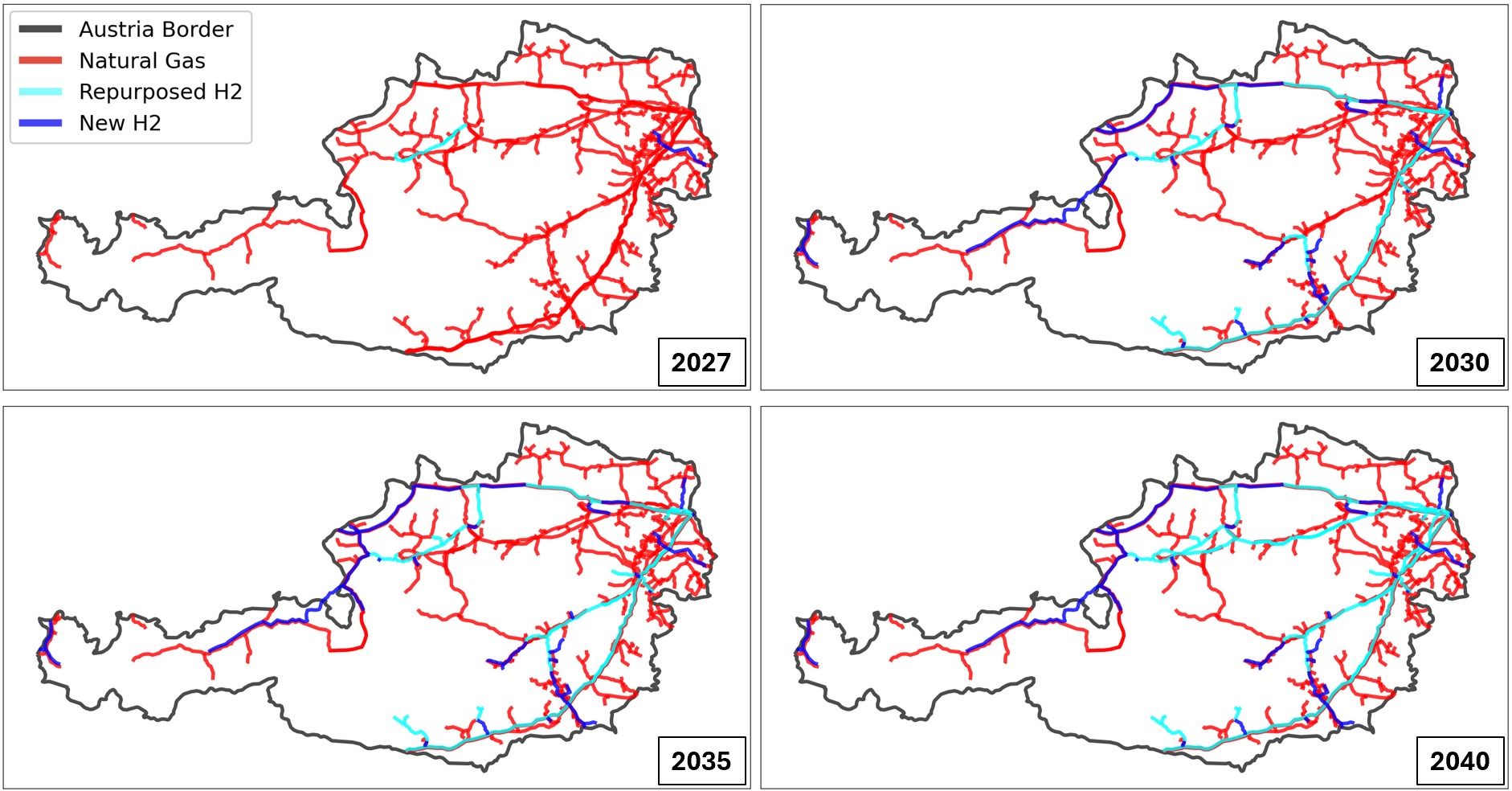}
    \caption{Austria's transition to an integrated hydrogen and natural gas infrastructure.}
    \label{fig:timeline}
\end{figure}

To illustrate the temporal development, Fig.~\ref{fig:timeline} shows a collage of the pipeline network evolution for 2027, 2030, 2035, and 2040. This visualization highlights the gradual expansion of hydrogen infrastructure and the repurposing of existing natural gas pipelines, providing an intuitive overview of the spatial and temporal dynamics of the resulting dataset.

This paper presents a methodology for developing a comprehensive, spatially explicit, and modeling-ready dataset of Austria's natural gas and hydrogen infrastructure. By integrating transmission and distribution networks with hydrogen expansion plans, the dataset enables consistent national and regional energy system analyses under network constraints.

By relying exclusively on publicly available data sources and reproducible processing steps, the dataset provides a transparent and updateable foundation for future infrastructure assessments. Due to limited public availability of detailed technical information, parts of the infrastructure, particularly at the distribution level, are represented using carefully approximated parameters. These approximations mainly affect local flow constraints and attribute resolution rather than the overall network topology. The dataset is therefore intended as a robust and flexible baseline representation rather than a finalized or static infrastructure model. Its modular, topology-consistent structure supports scenario-specific modifications, alternative hydrogen expansion and repurposing strategies, and future extensions to larger geographic scopes, including European-scale applications.

Overall, the dataset offers a high-resolution, reproducible foundation for transparent modeling of Austria's natural gas and hydrogen networks, supporting analyses of resilient, climate-neutral energy systems. The forthcoming release of the open-source QGas tool will enable further manipulation, extension, and integration of additional infrastructure components, making the dataset usable for a wide range of case studies such as strategic electrolyzer placement or analyses on new hydrogen corridors. To ensure transparency and reproducibility, the dataset presented in this paper is made publicly available on Zenodo~\cite{zenodo_dataset}.

\section*{Acknowledgment}

This work is part of the project iKlimEt (FO999910627), which has received funding in the framework of ``Energieforschung'', a research and technology program of the Klima- und Energiefonds.



\begin{thebibliography}{20}

\bibitem{kharel2018}
S.~Kharel and B.~Shabani,
``Hydrogen as a long-term large-scale energy storage solution to support renewables,''
\emph{Energies}, vol.~11, no.~10, p.~2825, 2018.

\bibitem{zwickl2024}
S.~Zwickl-Bernhard \emph{et al.},
``Modeling insights from the Austrian national gas grid under declining natural gas demand and increasing domestic renewable gas generation by 2040,''
\emph{Energy Reports}, vol.~11, pp.~1302--1317, 2024.

\bibitem{aryanpur2021}
V.~Aryanpur, B.~O'Gallach\'oir, H.~Dai, W.~Chen, and J.~Glynn,
``A review of spatial resolution and regionalisation in national-scale energy systems optimisation models,''
\emph{Energy Strategy Reviews}, vol.~37, 2021.

\bibitem{neumann2023}
F.~Neumann, E.~Zeyen, M.~Victoria, and T.~Brown,
``The potential role of a hydrogen network in Europe,''
\emph{Joule}, vol.~7, no.~8, pp.~1793--1817, 2023.

\bibitem{reuss2019}
M.~Reu{\ss} \emph{et al.},
``Modeling hydrogen networks for future energy systems: A comparison of linear and nonlinear approaches,''
\emph{International Journal of Hydrogen Energy}, vol.~44, no.~57, pp.~32136--32150, 2019.

\bibitem{entsog}
ENTSOG,
``Transparency Platform.''
[Online]. Available: \url{https://transparency.entsog.eu}.

\bibitem{pluta2022}
A.~Pluta \emph{et al.},
``SciGRID\_gas - data model of the European gas transport network,''
in \emph{2022 Open Source Modelling and Simulation of Energy Systems (OSMSES)}, 2022, pp.~1--7.

\bibitem{tag2025}
TAG GmbH,
``Technische Daten des TAG-Pipelinesystems,'' 2025.
[Online]. Available: \url{https://www.taggmbh.at/en/transmission-system/}.

\bibitem{gasconnect2025}
Gas Connect Austria,
``West-Austria Pipeline (WAG) -- technical and network information,'' 2025.
[Online]. Available: \url{https://www.gasconnect.at/en/network-information/our-network-in-detail/west-austria-pipeline}.

\bibitem{aggm_h2_roadmap}
AGGM Austrian Gas Grid Management AG,
``H2-Roadmap,'' 2024.
[Online]. Available: \url{https://www.aggm.at/en/energy-transition/h2-roadmap/}.

\bibitem{aggm_knep}
AGGM Austrian Gas Grid Management AG,
``Koordinierter Netzentwicklungsplan 2024 f\"ur die Gas-Fernleitungsinfrastruktur in \"Osterreich f\"ur den Zeitraum 2025--2034,'' 2024.

\bibitem{zwickl2023zenodo}
S.~Zwickl-Bernhard \emph{et al.},
``Replication package for the paper `Modeling insights from the Austrian national gas grid under declining natural gas demand and increasing domestic renewable gas generation by 2040',''
Zenodo, 2023.
[Online]. Available: \url{https://doi.org/10.5281/zenodo.10454605}.

\bibitem{quantschnig2025thesis}
M.~Quantschnig,
``Open-source gas infrastructure data for integrated energy system optimization models,''
Master's Thesis, Graz University of Technology, 2025.

\bibitem{osm}
OpenStreetMap Foundation,
``OpenStreetMap,'' 2025.
[Online]. Available: \url{https://www.openstreetmap.org}.

\bibitem{gem}
Global Energy Monitor,
``Global Gas Infrastructure Tracker,'' 2024.
[Online]. Available: \url{https://globalenergymonitor.org/projects/global-gas-infrastructure-tracker/}.

\bibitem{aggm_infra}
AGGM Austrian Gas Grid Management AG,
``Gas Infrastructure in Austria,'' 2025.
[Online]. Available: \url{https://www.aggm.at/en/gasgrid/infrastructure/}.

\bibitem{econtrol}
E-Control,
``Das \"osterreichische Gasnetz,'' 2025.
[Online]. Available: \url{https://www.e-control.at/konsumenten/das-gasnetz}.

\bibitem{dvgw2021}
DVGW -- Deutscher Verein des Gas- und Wasserfaches e.V.,
``Bestands- und Ereignisdatenerfassung Gas -- Ergebnisse aus den Jahren 2011 bis 2020,'' 2021.

\bibitem{gwg2011}
Republik \"Osterreich,
``Gaswirtschaftsgesetz 2011 (GWG 2011), BGBl. I Nr. 107/2011 idF BGBl. II Nr. 474/2012,'' 2012.

\bibitem{zenodo_dataset}
M.~Quantschnig, Y.~Werner, T.~Klatzer, and S.~Wogrin,
``Integrated natural gas and hydrogen dataset for energy system optimization models - Austria (2025 - 2040) (v1.0),''
Zenodo, EnInnov 2026, Graz, doi: 10.5281/zenodo.18411615.

\end{thebibliography}
\end{document}